\documentclass[11pt]{article}
\usepackage{amsmath}
\usepackage{amsfonts}
\usepackage{amssymb}
\usepackage{enumerate}
\usepackage{amsthm}
\usepackage{graphicx}
\usepackage{verbatim}
\newtheorem{myproposition}{Proposition}[section]
\newtheorem{mytheorem}[myproposition]{Theorem}
\newtheorem{mylemma}[myproposition]{Lemma}

\newtheorem{mycorollary}[myproposition]{Corollary}

\newtheorem{mydefinition}[myproposition]{Definition}

\newtheorem{myobservation}[myproposition]{Observation}

\newtheorem{myproblem}[myproposition]{Problem}

\def\ni{\noindent}

\def\gA{\mathcal{A}}
\def\zet{\mathbb{Z}}

\makeatletter
\def\imod#1{\allowbreak\mkern10mu({\operator@font mod}\,\,#1)}
\makeatother

\makeatletter
\def\jmod#1{\allowbreak\mkern10mu{\operator@font mod}\,\,#1}
\makeatother

\begin {document}
\title{Note on group distance magic graphs $G\times C_n$}

\author{Sylwia Cichacz\thanks{The author was supported by National Science Centre grant nr 2011/01/D/ST/04104.}\\
{\small Faculty of Applied Mathematics}\\
{\small AGH University of Science and Technology}\\
{\small Al. Mickiewicza 30, 30-059 Krak\'ow, Poland}}

\maketitle

\begin{abstract}
A  $\Gamma$-distance magic labeling
of a graph $G=(V,E)$ with $|V | = n$ is a bijection $f$ from $V$ to an Abelian group $\Gamma$ of
order $n$ such that the weight $w(x)=\sum_{y\in N_G(x)}f(y)$ of every vertex $x \in V$ is equal to the same element
$\mu \in \Gamma$, called the \emph{magic constant}.

In this paper we will show that if $G$ is a graph of order $n=2^{p}(2k+1)$ for some natural numbers $p$, $k$ such that  $\deg(v)\equiv c \imod {2^{p+2}}$ for some constant $c$ for any $v\in V(G)$, then there exists a $\Gamma$-distance magic labeling
for any Abelian group $\Gamma$ of order $4n$ for the direct product $G\times C_4$. Moreover if $c$ is even then there exists a $\Gamma$-distance magic labeling
for any Abelian group $\Gamma$ of order $8n$ for the direct product $G\times C_8$.

\end{abstract}
\noindent\textbf{Keywords:}

\section{Introduction and preliminaries}

All graphs considered in this paper are simple finite graphs. We use $V(G)$
for the vertex set and $E(G)$ for the edge set of a graph $G$. The \emph{\
neighborhood} $N(x)$ or more precisely $N_{G}(x)$, when needed, of a vertex $%
x$ is the set of vertices adjacent to $x$, and the \emph{degree} $d(x)$ of $%
x $ is $|N(x)|$, the size of the neighborhood of $x$. By $C_{n}$ we denote a
cycle on $n$ vertices.

A \emph{distance magic labeling} (also called \emph{sigma
labeling}) of a graph $G=(V,E)$ of order $n$ is a bijection $l
\colon V \rightarrow \{1, 2,\ldots  , n\}$ with the property that
there is a positive integer $\mu$ (called the \emph{magic constant}) such
that
$\sum_{y\in N_G(x)}l(y) =\mu$ for every $x \in V$. If a graph $G$ admits a distance magic labeling, then we say that $G$ is a  \emph{distance magic graph} (\cite{AFK}). The sum $\sum_{y\in N_G(x)}l(y)$ is called the \emph{weight of the vertex} $x$ and denoted by $w(x)$.\\

The concept of distance magic labeling has been motivated by the construction of magic squares. It is worth mentioning that finding an $r$-regular
distance magic labeling turns out equivalent to finding equalized
incomplete tournament $\mathrm{EIT}(n, r)$ \cite{FKK1}. In an \emph{equalized incomplete tournament} $\mathrm{EIT}(n, r)$ of $n$ teams with $r$
rounds,  every team plays exactly $r$ other teams and the total
strength of the opponents that team $i$ plays is $k$. Thus, it is easy to notice that finding an  $\mathrm{EIT}(n, r)$ is the same as finding a distance magic labeling of any $r$-regular graph on $n$ vertices. For a survey, we refer the reader to \cite{AFK}.

The following observations were independently proved:

\begin{myobservation}[\cite{Ji,MRS,Rao,Vi}] Let $G$ be an $r$-regular distance magic graph on $n$
vertices. Then $\mu = \frac{r(n+1)}{2}$.\end{myobservation}

\begin{myobservation}[\cite{Ji,MRS,Rao,Vi}] No $r$-regular graph with $r$ odd can be a distance
magic graph.\label{nieparzyste}\end{myobservation}

We recall two out of four standard graph products (see \cite{HIK}). Both, the \emph{lexicographic product} $G\circ H$ and
the \emph{direct product} $G\times H$ are graphs with the vertex set $V(G)\times
V(H)$. Two vertices $(g,h)$ and $(g^{\prime },h^{\prime })$ are adjacent in:

\begin{itemize}
\item $G\circ H$ if and only if either $g$ is adjacent with $g^{\prime }$ in $G$
or $g=g^{\prime }$ and $h$ is adjacent with $h^{\prime }$ in $H$;

\item $G\times H$ if $g$ is adjacent to $g^{\prime }$ in $G$ and $h$ is
adjacent to $h^{\prime }$ in $H$.
\end{itemize}
The graph $G\circ H$ is  also called the  \emph{composition} and denoted by $G[H]$ (see \cite{Har}). The  product $G\times H$  is also known as \emph{Kronecker
product}, \emph{tensor product}, \emph{categorical product} and \emph{graph conjunction} is the
most natural graph product. The direct product is commutative, associative, and it has several applications, for instance it may be used
as a model for concurrency in multiprocessor systems \cite{LaBa}. Some other
applications can be found in \cite{JKZ}.

Some graphs which are distance magic among
(some) products can be seen in \cite{ACG,ACPT,ACPT2,Be,Cic,Fro,MRS,RSP}.

\begin{mytheorem}[\cite{MRS}]  Let $r \geq 1$, $n \geq
3$, $G$ be an $r$-regular graph and $C_n$ be the cycle of length $n$.
The graph $G\circ C_n$ admits a distance magic labeling if and only if $n =
4$.
\end{mytheorem}

\begin{mytheorem}[\cite{MRS}]  Let $G$ be an arbitrary regular graph. Then $G\circ\overline{K}_{n}$ is distance
magic for any even $n$.
\end{mytheorem}

\begin{mytheorem}[\cite{ACPT}]  Let $G$ be an arbitrary regular graph. Then $G\times C_4$ is distance
magic.\label{C4A}
\end{mytheorem}

The following problem was posted  in \cite{AFK}.

\begin{myproblem}[\cite{AFK}] If $G$ is a non-regular graph, determine if there is a distance magic labeling of
      $G\circ C_4$.
\end{myproblem}

The similar problem for the direct product was stated in~\cite{ACG}:
\begin{myproblem}[\cite{ACG}] If $G$ is a non-regular graph, determine if there is a distance magic labeling of
      $G\times C_4$.
\end{myproblem}

Moreover it was proved that:

\begin{mytheorem}[\cite{ACG}]\label{ACG}
Let $m$ and $n$ be two positive integers such that $m\leq n$. The graph $K_{m,n}\times C_4$ is  a distance magic graph if and only if the following conditions hold:
\begin{enumerate}
  \item $m+n \equiv 0 \pmod 2$ and
  \item $m\geq \frac{\sqrt{2(8n+1)^2-1}-1}{8} -n.$
\end{enumerate}
\end{mytheorem}

Froncek in \cite{Fro} defined the notion of \emph{group distance magic graphs}, i.e. the graphs allowing the bijective labeling of vertices with elements of an Abelian group resulting in constant sums of neighbor labels.

\begin{mydefinition}
A $\Gamma$-distance magic labeling
of a graph $G=(V,E)$ with $|V | = n$ is a bijection $f$ from $V$ to an Abelian group $\Gamma$ of
order $n$ such that the weight $w(x)=\sum_{y\in N_G(x)}f(y)$ of every vertex $x \in V$ is equal to the same element
$\mu \in \Gamma$, called the magic constant.  A graph $G$ is called a \emph{group
distance magic graph} if there exists a $\Gamma $-distance magic labeling
for every Abelian group $\Gamma$ of order $|V(G)|$.

\end{mydefinition}

The connection between distance magic graphs and $\Gamma $-distance magic graphs is as follows. Let $G$ be a distance magic graph of order $n$ with the magic constant $\mu^{\prime }$. If we replace the label $n$ in a distance magic labeling for the graph $G$ by the label $0$, then we obtain a $\zet_n$-distance magic labeling for the graph $G$  with the magic constant $\mu\equiv \mu^{\prime } \imod n$. Hence every distance magic graph with $n$ vertices admits a $\mathbb{Z}_n$-distance magic labeling. Although a $\zet_n$-distance magic graph on $n$ vertices is
not necessarily a distance magic graph (see \cite{Fro}), it was proved that Observation~\ref{nieparzyste}  also holds for a $\zet_n$-distance magic labeling (\cite{CicFro}). 

\begin{myobservation}[\cite{CicFro}] Let $r$ be a positive odd integer. No $r$-regular graph  on $n$ vertices can be a $\zet_n$-distance
magic graph.
\end{myobservation}
 The following theorem was proved in \cite{Fro}:

\begin{mytheorem}[\protect\cite{Fro}]
\label{cart_cycle2}The Cartesian product $C_{m}\square C_{k}$, $m,k\geq 3$,
is a $\mathop\mathbb{Z}\nolimits_{mk}$-distance magic graph if and only if $km$
is even.
\end{mytheorem}

Froncek also showed that the graph $C_{2k}\square C_{2k}$ has a $(\mathop
\mathbb{Z}_2)^{2k}$-distance magic labeling for $k\geq 2$ and $\mu
=(0,0,\ldots ,0$) (\cite{Fro}).

Cichacz proved:
\begin{mytheorem}[\cite{Cic}]\label{lexmain}
Let $G$ be a graph of order $n$ and $\Gamma$  be an Abelian group of order $4n$. If $n=2^{p}(2k+1)$ for some natural numbers $p$, $k$ and   $\deg(v)\equiv c \imod {2^{p+1}}$ for some constant $c$ for any $v\in V(G)$, then there exists a $\Gamma$-distance magic labeling for the graph $G\circ C_4$.
\end{mytheorem}

It seems that the direct product is the natural choice among (standard)
products to deal with $\Gamma$-distance magic graphs and group distance
magic graphs in general. The reason for this is that the direct product is
suitable product if we observe graphs as categories. Hence it should perform
well with the product of (Abelian) groups, what the below theorem illustrates.
\begin{mytheorem}[\cite{ACPT2}]
If an $r_{1}$-regular graph $G_{1}$ is $\Gamma _{1}$-distance magic and an $%
r_{2}$ -regular graph $G_{2}$ is $\Gamma _{2}$-distance magic, then the
direct product $G_{1}\times G_{2}$ is $\Gamma _{1}\times \Gamma _{2}$%
-distance magic.\label{gr_pr}
\end{mytheorem}

However, dealing with this product is also most difficult in many
respects among standard products. For instance, $G\times H$ does not need to be
connected, even if both factors are. More precisely, $G\times H$ is
connected if and only if both $G$ and $H$ are connected and at least one of them is
non-bipartite \cite{Weich}. The main open problem concerning the direct product is the famous Hedetniemi's conjecture.
Hedetniemi  conjectured that for all graphs $G$ and $H$, $\chi (G\times H)=\min \{\chi (G),\chi (H)\}$, \cite{Hed}.\\
Anholcer at al. proved the following theorems:

\begin{mytheorem}[\cite{ACPT2}]
\label{C_4-all}If $G$ is an $r$-regular graph of order $n$%
, then lexicographic product $G\circ C_{4}$ is a group distance magic graph.
\end{mytheorem}

\begin{mytheorem}[\cite{ACPT2}]
\label{C_4-all2}If $G$ is an $r$-regular graph of order $n$%
, then direct product $G\times C_{4}$ is a group distance magic graph.
\end{mytheorem}

\begin{mytheorem}[\cite{ACPT2}]
\label{C_8-all}If $G$ is an $r$-regular graph of order $n$ for some even $r$%
, then direct product $G\times C_8$ is a group distance magic graph.
\end{mytheorem}

In this paper we prove the analogous theorems to Theorem~\ref{lexmain} for direct product $G\times C_{n}$ for $n\in\{4,8\}$.


\section{Direct product $G\times C_4$}
We start with the following lemma.

\begin{mylemma}\label{inne}
Let $G$ be a graph of order $n$ and $\Gamma$ be an arbitrary Abelian group of order $4n$ such that $\Gamma \cong \zet_{2^{p}}  \times \gA$ for $p\geq 1$ and some Abelian group $\gA$ of order $\frac{n}{2^{p-2}}$. If $\deg(v)\equiv c \imod {2^p}$ for some constant $c$ and any $v\in V(G)$, then there exists a $\Gamma$-distance magic labeling for the graph $G\times C_4$.
\end{mylemma}

\ni\textbf{Proof.}
 Let $V(G)=\{x_{0},x_{1},\ldots ,x_{n-1}\}$ be the
vertex set of $G$, let $C_{4}=u_{0}u_{1}u_2u_3u_{0}$, and $%
H=G\times C_4$. Let $v_{i,j}=(x_i,u_j)$ for $i=0,1,\ldots,n-1$ and  $j= 0,1,2,3$. Notice that if $x_{p}x_{q}\in E(G)$, then $v_{q,j}\in
N_{H}(v_{p,j})$ if and only if $j\in \{i-1,i+1\}$ (where the sum on the
second suffix is taken modulo $4$).

Using the isomorphism $\phi\colon \Gamma \rightarrow \gA\times \zet_{2^p}$, we identify every $g \in \Gamma$ with
its image $\phi(g) = (a_i,w)$, where $a_i \in \gA$ and $w \in \mathbb{Z}_{2^p}$, $i=0,1,\ldots,\frac{n}{2^{p-2}}-1$.\\

 Label the vertices of $H$ in the following way:

$$f(v_{i,j})=\left\{\begin{array}{lcc}
             \left(a_{\lfloor i \cdot 2^{-p+2}\rfloor},(2i+j)\jmod {2^{p-1}}\right) & \mathrm{for} & j=0,1, \\
             \left(0,2^p-1\right)-f(v_{i,j-2}) & \mathrm{for} & j=2,3 \\
          \end{array}\right.$$
for $i=0,1,\ldots,n-1$ and  $j= 0,1,2,3$. \\

 Notice that for every $i$
$$f(v_{i,0}) +f(v_{i,2})=f(v_{i,1}) + f(v_{i,3})=(0,2^{p}-1).$$
Since $\deg(v)\equiv c \imod {2^{p}}$ for any $v\in V(G) $, therefore the weight  of every $x\in V(H)$ is $w(x)=(0,-c)$.~\qed

\begin{mylemma}\label{lemma22}
Let $G$ be a graph of order $n$ and $\Gamma$ be an arbitrary Abelian group of order $4n$ such that $\Gamma \cong \zet_2 \times\zet_2 \times \gA$ for some Abelian group $\gA$ of order $n$.  If all vertices of $G$ have even degrees or all vertices of $G$ have odd degrees, then there exists a $\Gamma$-distance magic labeling for the graph $G\times C_4$.
\end{mylemma}

\ni\textbf{Proof.}
 Let $V(G)=\{x_{0},x_{1},\ldots ,x_{n-1}\}$ be the
vertex set of $G$, let $C_{4}=u_{0}u_{1}u_2 u_{3}u_{0}$, and $%
H=G\times C_4$. Let $v_{i,j}=(x_i,u_j)$ for $i=0,1,\ldots,n-1$ and  $j= 0,1,2,3$. Recall that if $x_{p}x_{q}\in E(G)$, then $v_{q,j}\in
N_{H}(v_{p,j})$ if and only if $j\in \{i-1,i+1\}$ (where the sum on the
second suffix is taken modulo $4$).  Since all vertices of $G$ have even degrees or all vertices of $G$ have odd degrees, thus $\deg(v)\equiv c \imod 2$ for some constant $c$ and any $v\in V(G)$

Using the isomorphism $\phi\colon \Gamma \rightarrow \gA\times \zet_{2}\times\zet_2$, we identify every $g \in \Gamma$ with
its image $\phi(g) = (a_i,j_1,j_2)$, where $j_1,j_2 \in \mathbb{Z}_2$ and $a_i \in \gA$, $i=0,1,\ldots,n-1$.\\

Label the vertices of $H$ in the following way:

$$f(v_{i,j})=\left\{\begin{array}{ccc}
            (a_i,0,0) & \mathrm{for} & j=0, \\
            (a_i,1,0) & \mathrm{for} & j=1, \\
            (-a_i,1,1) & \mathrm{for} & j=2, \\
            (-a_i,0,1) & \mathrm{for} & j=3 \\
          \end{array}\right.$$
for $i=0,1,\ldots,n-1$ and $j=0,1,2,3$.\\

Notice that for every $i=0,\ldots,n-1$
$$f(v_{i,0}) +f(v_{i,2})=f(v_{i,1}) + f(v_{i,3})=(0,1,1).$$
Therefore, for every $x\in V(H)$,
$$w(x) = (0,c,c).$$
\qed

\begin{mytheorem}\label{main}
Let $G$ be a graph of order $n$. If $n=2^{p}(2k+1)$ for some natural numbers $p$, $k$ and   $\deg(v)\equiv c \imod {2^{p+2}}$ for some constant $c$ for any $v\in V(G)$, then there exists a  group distance magic labeling for the graph $G\times C_4$.
\end{mytheorem}
\ni\textbf{Proof.}\\
The fundamental theorem of finite Abelian groups states that the finite Abelian group $\Gamma$ can be expressed as the direct sum of cyclic subgroups of prime-power order. This implies that $\Gamma \cong \mathbb{Z}_{2^{\alpha_0}} \times \mathbb{Z}_{p_{1}^{\alpha_1}}\times \mathbb{Z}_{p_{2}^{\alpha_2}}\times \ldots \times \mathbb{Z}_{p_{m}^{\alpha_m}}$ for some $\alpha_0>0$,
where $4n=2^{\alpha_0}\prod_{i=1}^m{p_i^{\alpha_i}}$ and $p_i$ for $i=1,\dots,m$ are not necessarily distinct primes. \\

\ni Suppose first that $\Gamma \cong \zet_2 \times \zet_2 \times \gA$ for some Abelian group $\gA$ of order $n$, then we are done by Lemma~\ref{lemma22}.
Observe now that the assumption  $\deg(v)\equiv c \imod {2^{p+2}}$  and the unique (additive) decomposition of any natural number $c$ into powers of $2$ imply that  there exist constants $c_1,c_2,\ldots,c_{p}$  such that
$\deg(v)\equiv c_i \imod {2^{i}}$,  for $i=1,2,\ldots,p+1$, for any $v\in V(G)$. Hence if $\Gamma \cong \mathbb{Z}_{2^{\alpha_0}} \times \gA$ for some $2 \leq \alpha_0\leq p+2$ and some Abelian group $\gA$ of order $\frac{4n}{2^{\alpha_0}}$, then
 we obtain by Lemma~\ref{inne} that there exists a $\Gamma$-distance magic labeling for the graph $G\times C_4$. \qed\\

The following observation shows that in the general case the condition on the degrees of the vertices $v \in V (G)$ in Theorem~\ref{main} is not necessary for the existence of a $\Gamma$-distance magic labeling of a graph $G\times C_4$:

 \begin{myobservation}\label{trojdzielne}
Let $G=K_{p,q,t}$ be a complete tripartite graph with all partite set odd, then $G\times C_4$ is a group distance magic graph.
\end{myobservation}
\ni\textbf{Proof.}\\
Since $n=p+q+t$ is odd $\Gamma\cong \gA\times \zet_2 \times\zet_2$ or $\Gamma\cong \gA\times \zet_4$ for some Abelian group $\gA$ of order $n=p+q+t$
    If  $\Gamma \cong \gA\times \zet_2 \times\zet_2$ for some Abelian group $\gA$ of order $n=p+q+t$, then there exists a $\Gamma$-distance magic labeling for the graph $K_{p,q,t}\times C_4$ by Lemma~\ref{lemma22}.
    Suppose now that $\Gamma \cong  \gA\times \zet_4 $ for some Abelian group $\gA$ of order $p+q+t$.
Let $K_{p,q,t}$ have the partition vertex sets $A=\{x_0, x_1,\ldots,x_{p-1}\}$, $B=\{y_0, y_1,\ldots,y_{q-1}\}$ and $C=\{u_0, u_1,\ldots,u_{t-1}\}$ and let $C_4=v_0v_1v_2v_3v_0$. Without loosing generality we can assume that $p\equiv q \pmod 4$ (by Pigeonhole Principle).
Using the isomorphism $\phi\colon \Gamma \rightarrow \gA\times \zet_{4}$, we identify every $g \in \Gamma$ with
its image $\phi(g) = (a_i,j)$, where $j \in \mathbb{Z}_4$ and $a_i \in \gA$, $i=0,1,\ldots,p+q+t-1$.

 Label the vertices of $K_{p,q,t}\times C_4$ in the following way:

$$f(x_i,v_j)=\left\{\begin{array}{ccc}
            (a_i,2j) & \mathrm{for} & j=0,1, \\
            \left(0,1\right)-f(x_i,v_{j-2}) & \mathrm{for} & j=2,3 \\
          \end{array}\right.$$
for $i=0,1,\ldots,p-1$ and $j=0,1,2,3$.\\

$$f(y_i,v_j)=\left\{\begin{array}{ccc}
            (a_{p+i},2j) & \mathrm{for} & j=0,1, \\
            \left(0,1\right)-f(y_i,v_{j-2}) & \mathrm{for} & j=2,3 \\
          \end{array}\right.$$
for $i=0,1,\ldots,q-1$ and $j=0,1,2,3$.\\

$$f(u_i,v_j)=\left\{\begin{array}{ccccl}
            (a_{p+q+i},2j) & \mathrm{for} & j=0,1, \\
            \left(0,1\right)-f(u_i,v_{j-2}) & \mathrm{for} & j=2,3, &\mathrm{if} & t\equiv p \pmod 4 \\
             \left(0,3\right)-f(u_i,v_{j-2}) & \mathrm{for} & j=2,3, &\mathrm{if} & t+2\equiv p \pmod 4 \\
          \end{array}\right.$$
for $i=0,1,\ldots,t-1$ and $j=0,1,2,3$.\\

Notice that
$$\begin{array}{c}
            f(x_i,v_j) +f(x_i,v_{j+2})=f(y_l,v_j) +f(y_l,v_{j+2})=(0,1),\\
          \end{array}
$$
for every $i=0,\ldots,p-1$, $l=0,\ldots,q-1$ and for $j=0,1$.\\
Whereas:
$$f(u_i,v_j) +f(u_i,v_{j+2})=\left\{
                               \begin{array}{lll}
                                 (0,1), & \hbox{if}& t\equiv p \pmod 4, \\
                                  (0,3), & \hbox{if}& t+2\equiv p \pmod 4, \\
                               \end{array}
                             \right.$$
for every $i=0,\ldots,t-1$ and for $j=0,1$.
Since $t\in\{1,3\}$ we obtain that $w(x)=(0,2)$  for every $x\in V(K_{p,q,t}\times C_4)$. \qed\\

Notice that a graph $K_{1,1+8\alpha}\times C_4$ is group distance magic for any $\alpha \in \mathbb{N}$ by Theorem~\ref{main}, although  a graph $K_{1,1+8\alpha}\times C_4$ is distance magic if and only if $\alpha=0$ by Theorem~\ref{ACG}.\\

It is worthy to mention, that it was proved that if $m$ is odd and $n$ is even, then a graph $K_{m,n}\circ C_4$ is group distance magic (see \cite{Cic}), however it is not longer true in the case of direct product,  what shows the bellow lemma. Recall that for a group $\Gamma$ an \emph{involution} $\iota\in \Gamma$ is a such element that $\iota\neq 0$ and $2\iota=0$.

 \begin{mylemma}\label{dwudzielne2}
Let $m$ and $n$ be two positive integers such that $m$ is odd and $n$ is even, the graph $K_{m,n}\times C_4$ is not a $\Gamma$-distance magic graph for any group $\Gamma$ of order $4m+4n$ having exactly one involution $\iota$.
\end{mylemma}
\ni\textbf{Proof.} Since there exists exactly one involution $\iota \in \Gamma$ notice that $\Gamma= \mathbb{Z}_{4}\times \gA$ for some group $\gA$ of order $m+n$. Let $K_{m,n}$ have the partition vertex sets $A=\{x_0, x_1,\ldots,x_{m-1}\}$ and $B=\{y_0, y_1,\ldots,y_{n-1}\}$ and let $C_4=v_0v_1v_2v_3v_0$.  Suppose that $\ell$ is a $\Gamma$-distance magic labeling of the graph $K_{m,n}\times C_4$ and
$\mu = w(x)$, for all vertices $x\in V(K_{m,n}\times C_4)$. Notice that $K_{m,n}\times C_4\cong 2K_{2m,2n}$. We can assume that $(x_i,v_j),(x_i,v_{j+2}),(y_l,v_{j+1}),(y_l,v_{j+3})\in V(K_{2p,2q}^j)$ for $i=0,1,\ldots,p-1$, $l=0,1,\ldots,q-1$ and $j=0,1$.  It is easy to observe that:
$$\mu=\sum_{i=0}^{m-1}\left(\ell(x_i,v_0)+\ell(x_i,v_2)\right)=\sum_{i=0}^{m-1}\left(\ell(x_i,v_1)+\ell(x_i,v_3)\right)=$$
$$=\sum_{l=0}^{n-1}\left(\ell(y_l,v_0)+\ell(y_l,v_2)\right)=\sum_{l=0}^{n-1}\left(\ell(y_l,v_1)+\ell(y_l,v_3)\right).$$

Thus $4\mu=\sum_{g\in\Gamma}g=\iota$. Since $m+n$ is odd and $\Gamma= \mathbb{Z}_{4}\times \gA$ such an element $\mu \in \Gamma$ does not exist, a contradiction.~\qed

\begin{mycorollary}\label{dwudzielne3}
Let $m$ and $n$ be two positive integers such that $m$ is odd and $n$ is even, the graph $K_{m,n}\times C_4$ is not a $\mathbb{Z}_{4m+4n}$-distance magic graph.
\end{mycorollary}
\ni\textbf{Proof.} Since there exists exactly one involution $\iota=2n+2m$ in $\mathbb{Z}_{4m+4n}$ we are done by Observation~\ref{dwudzielne2}.~\qed\\

We finish this section with the following theorem.

 \begin{mylemma}\label{dwudzielne2b}
Let $m$ and $n$ be two positive integers such that $m$ is odd and $n$ is even, the graph $K_{m,n}\times C_4$ is a $\Gamma$-distance  if and only if $\Gamma\cong \mathbb{Z}_{2}\times\mathbb{Z}_{2}\times \gA$ for a group $\gA$ of order $m+n$.
\end{mylemma}
\ni\textbf{Proof.} Since $|V(K_{m,n}\times C_4)|=4(m+n)$ and $m$ is even and $n$ is odd $\Gamma\cong \mathbb{Z}_{2}\times\mathbb{Z}_{2}\times \gA$ or $\Gamma= \mathbb{Z}_{4}\times \gA$ for some group $\gA$ of order $m+n$. If $\Gamma= \mathbb{Z}_{4}\times \gA$ then there does not exist a $\Gamma$-labeling of $K_{m,n}\times C_4$ by Lemma~\ref{dwudzielne2}. Suppose now that $\Gamma\cong \mathbb{Z}_{2}\times\mathbb{Z}_{2}\times \gA\cong \gA\times \mathbb{Z}_{2}\times\mathbb{Z}_{2}$ for a group $\gA$ of order $m+n$. Let $K_{m,n}$ have the partition vertex sets $A=\{x_0, x_1,\ldots,x_{m-1}\}$ and $B=\{y_0, y_1,\ldots,y_{n-1}\}$ and let $C_4=v_0v_1v_2v_3v_0$.  Using the isomorphism $\phi\colon \Gamma \rightarrow \gA\times \zet_{2}\times \zet_2$, we identify every $g \in \Gamma$ with
its image $\phi(g) = (a_i,j_1,j_2)$, where $j_1,j_2 \in \mathbb{Z}_4$ and $a_i \in \gA$, $i=0,1,\ldots,m+n-1$. Since $m+n$ is odd without loosing the generality that $a_0=0$ and $a_{n+1}=-a_{n}\neq 0$.\\

Label the vertices of $K_{m,n}\times C_4$ in the following way:\\
$f(x_0,v_0)=(0,0,0)$, $f(x_0,v_1)=(0,1,0)$, $f(x_0,v_2)=(0,0,1)$, $f(x_0,v_3)=(0,1,1)$,\\
$f(y_0,v_0)=(a_n,1,0)$, $f(y_0,v_1)=(a_n,0,0)$, $f(y_0,v_2)=(-a_n,1,0)$, $f(y_0,v_3)=(-a_n,1,1)$,\\
$f(y_1,v_0)=(a_{n+1},0,0)$, $f(y_1,v_1)=(a_{n+1},0,1)$, $f(y_1,v_2)=(-a_{n+1},0,1)$, $f(y_1,v_3)=(-a_{n+1},1,1)$.\\

$$f(x_i,v_j)=\left\{\begin{array}{ccc}
            (a_i,0,0) & \mathrm{for} & j=0, \\
            (a_i,1,0) & \mathrm{for} & j=1, \\
            (-a_i,1,1) & \mathrm{for} & j=2, \\
            (-a_i,0,1) & \mathrm{for} & j=3 \\
          \end{array}\right.$$
for $i=1,2,\ldots,m-1$ and $j=0,1,2,3$.\\

$$f(y_i,v_j)=\left\{\begin{array}{ccc}
           (a_{m+i},0,0) & \mathrm{for} & j=0, \\
            (a_{m+i},1,0) & \mathrm{for} & j=1, \\
            (-a_{m+i},1,1) & \mathrm{for} & j=2, \\
            (-a_{m+i},0,1) & \mathrm{for} & j=3 \\
          \end{array}\right.$$
for $i=2,3,\ldots,n-1$ and $j=0,1,2,3$.\\

Notice that
$$\begin{array}{c}
            f(x_i,v_j) +f(x_i,v_{j+2})=f(y_l,v_j) +f(y_l,v_{j+2})=(0,1,1),\\
          \end{array}
$$
for every $i=1,2,\ldots,m-1$, $l=2,3\ldots,n-1$ and for $j=0,1$.\\
We obtain that $w(x_i,v_j)=(n-2)(0,1,1)+(0,0,1)=(0,0,1)$ for $i=0,1,\ldots,m-1$ and $j=0,1,2,3$ and $w(y_l,v_j)=(m-1)(0,1,1)+(0,0,1)=(0,0,1)$ for $i=0,1,\ldots,m-1$ and $j=0,1,2,3$.~\qed\\


\section{Direct product $G\times C_8$}

In this section we   show that some  direct products $G\times C_8$ are group distance magic. Used constructions are similar to those by
Anholcer at al. in~\cite{ACPT2}.  We start with the following lemma:

\begin{mylemma}\label{lemma33}
Let $G$ be a graph of order $n$ and $\Gamma$ be an arbitrary Abelian group of order $4n$ such that $\Gamma \cong \zet_2 \times\zet_2 \times \gA$ for some Abelian group $\gA$ of order $2n$.  If all vertices of $G$ have even degrees, then there exists a $\Gamma$-distance magic labeling for the graph $G\times C_8$.
\end{mylemma}

\ni\textbf{Proof.}
 Let $V(G)=\{x_{0},x_{1},\ldots ,x_{n-1}\}$ be the
vertex set of $G$, let $C_{4}=u_{0}u_{1}\ldots u_{7}u_{0}$, and $%
H=G\times C_8$. Let $v_{i,j}=(x_i,u_j)$ for $i=0,1,\ldots,n-1$ and  $j= 0,1,\ldots,7$. Notice that if $x_{p}x_{q}\in E(G)$, then $v_{q,j}\in
N_{H}(v_{p,i})$ if and only if $j\in \{i-1,i+1\}$ (where the sum on the
first suffix is taken modulo $8$).  Since all vertices of $G$ have even degrees, thus $\deg(v)=2l_v$  any $v\in V(G)$.

Using the isomorphism $\phi\colon \Gamma \rightarrow \gA\times \zet_{2}\times\zet_2$, we identify every $g \in \Gamma$ with
its image $\phi(g) = (a_i,j_1,j_2)$, where $j_1,j_2 \in \mathbb{Z}_2$ and $a_i \in \gA$, $i=0,1,\ldots,2n-1$.\\

For $j\in \{0,1,\ldots ,n-1\}$ we set

\begin{equation*}
f (v_{i,j})=\left\{
\begin{array}{lcl}
(a_{2i+j},0,0), & \text{if} & i\in \{0,1\}, \\
(a_{2i+j-4},0,1), & \text{if} & i\in \{4,5\}, \\
(0,1,1)-f(v_{i,j-2}), & \text{if} & i\in \{2,3,6,7\}.%
\end{array}%
\right.
\end{equation*}%
Clearly, $f \colon V(C_{8}\times G)\rightarrow \Gamma $ is a bijection and $f(v_{i,j})+f(v_{i,j+2})=(0,y_{j})$, where $y_{j}\in \{(1,1),(1,0)\}$,
and so $2y_{j}=(0,0)$. Hence for every $i\in
\{0,1,\ldots ,n-1\}$ and $j\in \{0,1,\ldots ,7\}$  we get
\begin{eqnarray*}
w(v_{i,j}) &=&\sum_{x_{p}\in N_{G}(x_{i})}(f(v_{p,j-1})+f(v_{p,j+1}))=\sum_{x_{p}\in N_{G}(x_{i})}(0,y_{j-1})= \\
&=&l_{v_{i,j}}(0,0,0)=(0,0,0)
\end{eqnarray*}%
and $G\times C_8$ is $\Gamma $-distance magic.~\qed

\begin{mytheorem}\label{C_8}
Let $G$ be a graph of order $n$. If $n=2^{p}(2k+1)$ for some natural numbers $p$, $k$ and   $\deg(v)\equiv 2c \imod {2^{p+2}}$ for some constant $c$ for any $v\in V(G)$, then there exists a  group distance magic labeling for the graph $G\times C_8$.
\end{mytheorem}

\noindent \textit{Proof.}\ Let $V(G)=\{x_{0},x_{1},\ldots ,x_{n-1}\}$ be the
vertex set of $G$, let $C_{8}=u_{0}u_{1}\ldots u_{7}u_{0}$, and $%
H=G\times C_8$. Let $v_{i,j}=(x_i,u_j)$ for $i=0,1,\ldots,n-1$ and  $j= 0,1,\ldots,7$. Notice that if $x_{p}x_{q}\in E(G)$, then $v_{q,j}\in
N_{H}(v_{p,j})$ if and only if $j\in \{i-1,i+1\}$ (where the sum on the
first suffix is taken modulo $8$).

Recall that the assumption  $\deg(v)\equiv c \imod {2^{p+2}}$  and the unique (additive) decomposition of any natural number $c$ into powers of $2$ imply that  there exist constants $c_1,c_2,\ldots,c_{p}$  such that
$\deg(v)\equiv 2c_i \imod {2^{i}}$,  for $i=1,2,\ldots,p+1$, for any $v\in V(G)$.

We are going to consider three cases,
depending on the structure of $\Gamma$.\newline

\noindent \textit{Case 1:} $\Gamma\cong \mathop\mathbb{Z}\nolimits_{2}\times %
\mathop\mathbb{Z}\nolimits _{2}\times \mathcal{A}$ for some Abelian group of
order $2n$.\\
There exists a $\Gamma$-distance magic labeling of $G\times C_8$ by Lemma~\ref{lemma33}.\\

\noindent \textit{Case 2:} $\Gamma \cong \mathop\mathbb{Z}\nolimits
_{4}\times \mathcal{A}$ for some Abelian group $\mathcal{A}$ of order $2n$.

Using the isomorphism $\phi\colon \Gamma \rightarrow \gA\times \zet_{4}$, we identify every $g \in \Gamma$ with
its image $\phi(g) = ( a_i,w)$, where $w \in \mathbb{Z}_{4}$ and $a_i \in \gA$, $i=0,1,\ldots,2n-1$.\\

For $i\in \{0,1,\ldots ,n-1\}$ we define

\begin{equation*}
f(v_{i,j})=\left\{
\begin{array}{lcl}
(a_{2i+j},0), & \text{if} & j\in \{0,1\}, \\
(a_{2i+j-4},2), & \text{if} & j\in \{4,5\}, \\
(0,3)-f(v_{i,j-2}), & \text{if} & j\in \{2,3,6,7\}.%
\end{array}
\right.
\end{equation*}
Again $f\colon V(G\times C_{8})\rightarrow \Gamma $ is obviously a bijection and $%
f (v_{i,j})+f (v_{i,j+2})=(y_i,0)$, where $y_j\in\{1,3\}$, and thus $%
2y_j=2$. Since $\deg(v)\equiv 2c_2 \imod {2^{2}}$  for any $v\in V(G)$, for every $i\in \{0,1,\ldots
,n-1\}$ and $j\in \{0,1,\ldots ,7\}$ we get
\begin{equation*}
w(v_{i,j})=\sum_{x_{p}\in N_{G}(x_{i})}(f(v_{p,j-1})+f
(v_{p,j+1}))=\sum_{x_{p}\in N_{G}(x_{i})}(0,y_{j-1})=c_2%
(0,2)=(0,2c_2)
\end{equation*}
and $G\times C_{8}$ is $\Gamma $-distance magic.\newline

\noindent \textit{Case 3:} $\Gamma \cong \mathop\mathbb{Z}%
\nolimits_{2^{\alpha }}\times \mathcal{A}$ for $2<\alpha\leq p$ and some Abelian
group $\mathcal{A}$ of order $\frac{n}{2^{\alpha -3}}$.

Using the isomorphism $\phi\colon \Gamma \rightarrow \gA\times \zet_{2^{\alpha }}$, we identify every $g \in \Gamma$ with
its image $\phi(g) = (a_i,w)$, where $w \in \mathop\mathbb{Z}%
\nolimits_{2^{\alpha }}$ and $a_{i}\in \mathcal{A}$ for $i\in \{0,1,\ldots ,%
\frac{n}{2^{\alpha -3}}-1\}$. For $i\in \{0,1,\ldots ,\frac{n}{2^{\alpha -3}}%
-1\}$ define the following labeling $f$:\\

\begin{equation*}
f(v_{i,j})=\left\{
\begin{array}{lcl}
\left( a_{\lfloor i\cdot2^{-\alpha +3}\rfloor },(2i+j)(\mathop{\rm mod}\nolimits{2^{\alpha -2}})\right) , & \text{if} & j\in \{0,1\}, \\
\left(0, 2^{\alpha -1})+f(v_{i,j-4}\right) , & \text{if} & j\in
\{4,5\}, \\
\left( 0,2^{\alpha }-1\right) -f(v_{i,j-2}), & \text{if} & j\in
\{2,3,6,7\}.%
\end{array}%
\right.
\end{equation*}%
As in previous cases $f:V(G\times C_{8})\rightarrow \Gamma $ is a bijection
and $f(v_{i,j})+f(v_{i,j+2})=(0,y_{j})$ for some $y_{j}\in
\{2^{\alpha -1}-1,2^{\alpha }-1\}$. Thus $2(f (v_{i,j})+f(v_{i,j+2}))=(0,2y_{j})=(0,-2)$. For every $j\in \{0,1,\ldots ,7\}$
and $i\in \{0,1,\ldots ,n-1\}$ we get
\begin{eqnarray*}
w(v_{i,j}) &=&\sum_{x_{p}\in N_{G}(x_{i})}(f(v_{p,j-1})+f(v_{p,j+1}))=\sum_{x_{p}\in N_{G}(x_{i})}(0,y_{j-1})= \\
&=&c_{\alpha}(0,-2)=(0,-2c_\alpha)
\end{eqnarray*}%
and $G\times C_{8}$ is $\mathop\mathbb{Z}\nolimits_{2^{\alpha }}\times
\mathcal{A}$-distance magic since $r$ is even.~\hfill \rule{0.1in}{0.1in}%
\medskip

In the proof of the below observation we use similar methods to those presented in~\cite{Be}.
 \begin{myobservation}\label{dwudzielne3}
Let $m$ and $n$ be two positive integers such that $m\leq n$. If the graph $K_{m,n}\times C_8$ is  a distance magic graph, then the conditions hold:
\begin{itemize}
  \item $m+n \equiv 0 \pmod 2$ and
  \item $m\geq \frac{\sqrt{2(16n+1)^2-1}-1}{16} -n.$
\end{itemize}
\end{myobservation}
\ni\textbf{Proof.} Let $K_{m,n}$ have the partition vertex sets $A=\{x_0, x_1,\ldots,x_{m-1}\}$ and $B=\{y_0, y_1,\ldots,y_{n-1}\}$ and let $C_4=v_0v_1v_2\ldots v_7v_0$.  Suppose that $\ell$ is a distance magic labeling of the graph $K_{m,n}\times C_8$ and
$\mu = w(x)$, for all vertices $x\in V(K_{m,n}\times C_8)$.  We can assume that $(x_i,v_j),(x_i,v_{j+2}),(y_l,v_{j+1}),(y_l,v_{j+3})\in V(K_{2p,2q}^j)$ for $i=0,1,\ldots,p-1$, $l=0,1,\ldots,q-1$ and $i=0,1,\ldots,7$.  It is easy to observe that:
$$\mu=\sum_{i=0}^{m-1}\left(\ell(x_i,v_0)+\ell(x_i,v_2)\right)=\sum_{i=0}^{m-1}\left(\ell(x_i,v_1)+\ell(x_i,v_3)\right)=$$
$$=\sum_{i=0}^{m-1}\left(\ell(x_i,v_4)+\ell(x_i,v_6)\right)=\sum_{i=0}^{m-1}\left(\ell(x_i,v_5)+\ell(x_i,v_7)\right)=$$
$$=\sum_{l=0}^{n-1}\left(\ell(y_l,v_0)+\ell(y_l,v_2)\right)=\sum_{l=0}^{n-1}\left(\ell(y_l,v_1)+\ell(y_l,v_3)\right)=$$
$$=\sum_{l=0}^{n-1}\left(\ell(y_l,v_4)+\ell(y_l,v_6)\right)=\sum_{l=0}^{n-1}\left(\ell(y_l,v_5)+\ell(y_l,v_7)\right)=$$
$$=\sum_{x\in V(K_{m,n}\times C_8)}\frac{l(x)}{4}=\sum_ {i=1}^{8n+8m}\frac{i}{8}=\frac{(8m+8n)(8m+8n+1)}{16},$$
which implies that $m+n \equiv 0 \pmod 2$.\\
Notice that $\sum_{i=0}^{m-1}\left(\ell(x_i,v_0)+\ell(x_i,v_1)+\ldots+\ell(x_i,v_6)+\ell(x_i,v_7)\right)\leq \sum_{i=1}^{8m}(i+8n)=4m(8m+16n+1)$, thus $\mu\leq m(8m+16n+1)$. Which implies $(m+n)(8m+8n+1)\leq2m(8m+16n+1)$ and therefore:
$$2[2m+(2n+\frac{1}{8})]^2\geq n(8n+1)+ 2(2n+\frac{1}{8})^2=(4n+\frac{1}{4})^2-\frac{1}{16}.$$
That is:
$$1\geq (16n+1)^2-2(8m+8n+\frac{1}{2})^2,$$
Therefore, either $2=2(16n+1)^2-(16m+16n+1)^2$ or $m\geq \frac{\sqrt{2(16n+1)^2-1}-1}{16} -n.$\\
Suppose that $1\geq 2(16n+1)^2-(16m+16n+1)^2$, then the diophantine equation $2=2x^2-y^2$ needs to have a solution such that $x$ and $y$ are both odd, a contradiction.~\qed\\

Notice that a graph $K_{2,2+16\alpha}\times C_4$ is group distance magic for any $\alpha \in \mathbb{N}$ by Theorem~\ref{C_8}, although  a graph $K_{2,2+16\alpha}\times C_4$ is not distance magic for $\alpha>0$ by Obseravtion~\ref{dwudzielne3}.\\

\end{document}